\documentclass{amsart}

\usepackage{amscd,amsmath,xypic,amssymb,combelow,tikz-cd,etoolbox,calligra,mathrsfs,enumitem,mathtools,epigraph,hyperref}

\makeatletter
\patchcmd{\@settitle}{\uppercasenonmath\@title}{}{}{}
\makeatother

\xyoption{all}
\CompileMatrices

\makeatletter
\@addtoreset{equation}{section}
\makeatother

\newtheorem{theorem}[subsection]{Theorem}
\newtheorem{problem}[subsection]{Problem}

\newtheorem{proposition}[subsection]{Proposition}
\newtheorem{lemma}[subsection]{Lemma}

\newtheorem{definition}[subsection]{Definition}

\newtheorem{claim}[subsection]{Claim}

\newtheorem{remark}[subsection]{Remark}

\def\loccit{\emph{loc. cit. }}

\def\fgl{{\mathfrak{gl}}}

\def\BA{{\mathbb{A}}}
\def\BC{{\mathbb{C}}}

\def\BN{{\mathbb{N}}}

\def\BP{{\mathbb{P}}}

\def\BQ{{\mathbb{Q}}}
\def\BZ{{\mathbb{Z}}}

\def\CA{{\mathcal{A}}}
\def\CB{{\mathcal{B}}}

\def\CO{{\mathcal{O}}}

\def\CS{{\mathcal{S}}}

\def\CV{{\mathcal{V}}}

\def\ph{\varphi}

\def\and{\textrm{ }\&\textrm{ }}

\def\sym{\textrm{sym}}
\def\Sym{\textrm{Sym}}
\def\esym{\emph{sym}}
\def\eSym{\emph{Sym}}

\def\oP{\overline{P}}

\def\loc{\text{loc}}
\def\eloc{\emph{loc}}

\def\ocomm{\overline{\comm}}
\def\comm{\text{Comm}}
\def\nilp{\text{semi-nilp}}
\def\enilp{\emph{\nilp}}

\def\oP{\bar{P}}

\def\oH{\bar{H}}

\begin{document}

\title[An integral form of quantum toroidal $\fgl_1$]{\Large{\textbf{An integral form of quantum toroidal $\fgl_1$}}}

\author[Andrei Negu\cb t]{Andrei Negu\cb t}
\address{MIT, Department of Mathematics, Cambridge, MA, USA}
\address{Simion Stoilow Institute of Mathematics, Bucharest, Romania}
\email{andrei.negut@gmail.com}

\maketitle

\begin{abstract} We consider the (direct sum over all $n$ of the) $K$-theory of the semi-nilpotent commuting variety of $\fgl_n$, and describe its convolution algebra structure in two ways: the first as an explicit shuffle algebra (i.e. a particular $\BZ[q_1^{\pm 1}, q_2^{\pm 1}]$-submodule of the equivariant $K$-theory of a point) and the second as the $\BZ[q_1^{\pm 1}, q_2^{\pm 1}]$-algebra generated by certain elements $\{\oH_{n,d}\}_{(n,d) \in \BN \times \BZ}$.

\end{abstract}

$$$$

\section{Introduction}

\medskip

\subsection{} 

Moduli spaces of quiver representations and moduli spaces of sheaves are both important settings for geometric representation theory. Moreover, they are very closely connected, in that one can see the same phenomena occur for both classes of moduli spaces. Arguably nowhere is this more apparent than in the case of the Jordan quiver (namely the quiver with one vertex and one loop), which corresponds to sheaves on $\mathbb{A}^2$. To be more specific, consider the commuting stack
$$
\text{Comm}_n = \Big\{(X,Y) \in \text{Mat}_{n \times n}^{\times 2} \text{ s.t. } [X,Y] = 0\Big\} \Big / GL_n
$$
where the action of $GL_n$ is by simultaneous conjugation of the matrices $X,Y$. From the point of view of quivers, $\text{Comm}_n$ is the cotangent bundle of the stack $\text{Mat}_{n \times n}/GL_n$ of $n$-dimensional representations of the Jordan quiver. From the point of view of sheaves, a point of $\text{Comm}_n$ describes a length $n$ sheaf on $\BA^2$, as the commuting endomorphisms $X$ and $Y$ encode an action of $\CO_{\BA^2}$. Let us consider
\begin{equation}
\label{eqn:k intro}
K = \bigoplus_{n=0}^{\infty} K_{\BC^* \times \BC^*}(\comm_n)
\end{equation}
the (0-th) equivariant algebraic $K$-theory groups of all commuting stacks considered together. The torus $\BC^* \times \BC^*$ acts by rescaling the matrices $X$ and $Y$ independently, and thus $K$ is a $\BZ[q_1^{\pm 1}, q_2^{\pm 1}]$-module, where $q_1,q_2$ denote the standard characters of $\BC^* \times \BC^*$. As explained in \cite{SV Hilb}, there is a convolution algebra structure on $K$ which is additive in $n$ (we will not need to review the construction in the present paper, but the interested reader may find an overview in \cite[Section 2.3]{Quiver 1}).

\medskip

\subsection{} Upon localization with respect to the fraction field of $\BZ[q_1^{\pm 1}, q_2^{\pm 1}]$, the algebra
\begin{equation}
\label{eqn:loc intro}
K_{\loc} = K \bigotimes_{\BZ[q_1^{\pm 1}, q_2^{\pm 1}]} \BQ(q_1,q_2)
\end{equation}
is a well-known object in representation theory: it was shown in \cite{SV Hilb} to match the elliptic Hall algebra of \cite{BS}, in \cite{FT, S} to match half of quantum toroidal $\fgl_1$ (also known as the Ding-Iohara-Miki algebra), and in \cite{Shuf} to match the shuffle algebra of \cite{FHHSY}. However, if one has derived categories (or any other categorification) in mind, knowing $K_{\loc}$ is not good enough. Instead, one would hope to solve the following.

\medskip

\begin{problem}
\label{prob}

Describe $K$ as a $\BZ[q_1^{\pm 1}, q_2^{\pm 1}]$-algebra.

\end{problem}

\medskip

\noindent Although certain aspects of Problem \ref{prob} have been studied (\cite{W, Z2}), the bad news is that we do not yet know a complete solution. The good news is that in Theorem \ref{thm:main}, we provide a complete solution to a closely related problem, which is relevant to the setting of categorified knot invariants and affine Hecke algebras studied in \cite{GN, GNR}. To set up this closely related problem, let us note that the commuting stack has three variants of interest to us, namely
\begin{equation}
\label{eqn:inclusion intro}
\comm_n^{\text{nilp}} \subset \comm_n^{\nilp} \subset \comm_n
\end{equation}
where the stack on the left consists of pairs of nilpotent commuting matrices $(X,Y)$, while the stack in the middle allows $X$ to be arbitrary but requires $Y$ to be nilpotent. The (direct sums over all $n \in \BN$ of the) $\BC^* \times \BC^*$ equivariant algebraic $K$-theory groups of the stacks above will be denoted by
\begin{equation}
\label{eqn:three}
K^{\text{nilp}} \longrightarrow K^{\nilp} \longrightarrow K
\end{equation}
The maps above are simply the direct image maps induced by \eqref{eqn:inclusion intro}, and they are actually algebra homomorphisms with respect to the convolution product (indeed, $K^{\text{nilp}}$ and $K^{\nilp}$ are algebras by the exact same construction as $K$ of \eqref{eqn:k intro}). All three algebras in \eqref{eqn:three} have the same localization, i.e. are isomorphic upon tensoring with $\BQ(q_1,q_2)$, but the middle one will be described explicitly before localization. 

\medskip

\begin{theorem}
\label{thm:main}

We have an isomorphism $\iota^{\enilp} : K^{\enilp} \xrightarrow{\sim} \CS$, where
$$
\CS \subset \bigoplus_{n=0}^{\infty} \BZ[q_1^{\pm 1},q_2^{\pm 1}][z_1^{\pm 1}, \dots, z_n^{\pm 1}]^{\esym}
$$
is the $\BZ[q_1^{\pm 1},q_2^{\pm 1}]$-submodule determined by the conditions of Definition \ref{def:shuffle}, and made into an algebra via the shuffle product \eqref{eqn:shuf prod}. 

\end{theorem}

\medskip

\subsection{} Our starting point in the analysis of $K^{\nilp}$ is the fact (proved in \cite{SV gen}) that it is generated as an algebra by the $K$-theory groups of the closed substacks
\begin{equation}
\label{eqn:substack intro}
\Big(\text{Mat}_{n \times n} \times \{0\}\Big) \Big/ GL_n \subset \text{Comm}_n^{\nilp}
\end{equation}
as $n$ ranges over $\BN$. The isomorphism $\iota^{\nilp}$ of Theorem \ref{thm:main} maps the $K$-theory group of the substack \eqref{eqn:substack intro} to
\begin{equation}
\label{eqn:ideal intro}
\BZ[q_1^{\pm 1}, q_2^{\pm 1}][z_1^{\pm 1},\dots,z_n^{\pm 1}]^{\sym} \cdot F_n \subset \CS
\end{equation}
where
\begin{equation}
\label{eqn:fn}
F_n = \prod_{1\leq i, j \leq n} \left(1 - \frac {q_2z_i}{z_j} \right)
\end{equation}
The elements $F_n$ were first studied in \cite{FHHSY}, and we prove the surjectivity of the map $\iota^{\nilp}$ by showing that the elements of \eqref{eqn:ideal intro} also generate $\CS$, as $n$ ranges over $\BN$. As the injectivity of $\iota^{\nilp}$ was established in \cite{VV}, this proves Theorem \ref{thm:main}. 

\medskip

\noindent Inspired by the elliptic Hall algebra of \cite{BS}, it was shown in \cite{Shuf} that we have the following equality of $\BQ(q_1,q_2)$-vector spaces
\begin{equation}
\label{eqn:pbw loc intro}
\CS_{\loc} := \CS \bigotimes_{\BZ[q_1^{\pm 1}, q_2^{\pm 1}]} \BQ(q_1,q_2) = \bigoplus_{\frac {d_1}{n_1} \leq \dots \leq \frac {d_k}{n_k}} \BQ(q_1,q_2) \cdot \oH_{n_1,d_1} \dots \oH_{n_k,d_k}
\end{equation}
for certain elements $\{\oH_{n,d}\}_{(n,d) \in \BN \times \BZ}$ of $\CS_{\loc}$,  that we will recall in Section \ref{sec:pbw}. We will then prove the following stronger version of the decomposition \eqref{eqn:pbw loc intro}.

\medskip

\begin{lemma}
\label{lem:pbw}

We have the following equality of $\BZ[q_1^{\pm 1}, q_2^{\pm 1}]$-modules
\begin{equation}
\label{eqn:pbw}
\CS = \bigoplus_{\frac {d_1}{n_1} \leq \dots \leq \frac {d_k}{n_k}} \BZ[q_1^{\pm 1}, q_2^{\pm 1}] \cdot \oH_{n_1,d_1} \dots \oH_{n_k,d_k}
\end{equation}

\end{lemma}

\medskip

\subsection{}  I would like to thank Eugene Gorsky, Olivier Schiffmann and Alexander Tsymbaliuk for numerous wonderful conversations about this problem and its many facets. I gratefully acknowledge NSF grant DMS-$1845034$, as well as support from the Alfred P.\ Sloan Foundation and the MIT Research Support Committee.

\bigskip

\section{The (semi-nilpotent) $K$-theoretic Hall algebra}

\medskip

\subsection{} Let us consider the \textbf{commuting variety}
\begin{equation}
\label{eqn:inclusion}
\ocomm_n \stackrel{i_n}\hookrightarrow \BA^{2n^2}
\end{equation}
\footnote{Strictly speaking, one should think of $\ocomm_n$ as the derived subscheme of $\BA^{2n^2}$ cut out by the Koszul complex of the system of $n^2$ equations $[X,Y] = 0$, but we will not need this subtlety.} consisting of pairs of commuting $n \times n$ matrices $X,Y$. We consider the action
$$
G_n := \BC^* \times \BC^* \times GL_n \curvearrowright \BA^{2n^2}
$$
given by
$$
(t_1,t_2,g) \cdot (X,Y) = \left( \frac 1{t_1} gXg^{-1}, \frac 1{t_2} gYg^{-1} 
\right)
$$
which preserves $\ocomm_n$. Thus, \eqref{eqn:inclusion} induces a map on equivariant $K$-theory
\begin{equation}
\label{eqn:inclusion push}
K_{G_n} (\ocomm_n) \xrightarrow{i_{n*}} K_{G_n} (\BA^{2n^2})
\end{equation}
The \textbf{commuting stack} is $\comm_n = \ocomm_n/GL_n$, and its $K$-theory is given by
\begin{equation}
\label{eqn:k of stack}
K_{\BC^* \times \BC^*}(\comm_n) = K_{G_n}(\ocomm_n)
\end{equation}
which explains our interest in the map \eqref{eqn:inclusion push}.

\medskip

\subsection{} If we let $\circ \in \BA^{2n^2}$ denote the origin, then the following restriction map
\begin{equation}
\label{eqn:restriction}
K_{G_n} (\BA^{2n^2}) \stackrel{|_\circ}\cong K_{G_n} (\text{pt}) = \BZ[q_1^{\pm 1}, q_2^{\pm 2}][z_1^{\pm 1}, \dots, z_n^{\pm 1}]^{\sym}
\end{equation}
is an isomorphism, where $q_1,q_2$ denote the standard characters of $\BC^* \times \BC^*$, $z_1,\dots,z_n$ denote the standard characters of a maximal torus of $GL_n$, and ``sym" denotes symmetric Laurent polynomials in $z_1,\dots,z_n$. Composing \eqref{eqn:inclusion push} with \eqref{eqn:restriction} yields
$$
K_{G_n} (\ocomm_n) \xrightarrow{\iota_n} \BZ[q_1^{\pm 1}, q_2^{\pm 2}][z_1^{\pm 1}, \dots, z_n^{\pm 1}]^{\sym}
$$
as a map of $\BZ[q_1^{\pm 1}, q_2^{\pm 1}]$-modules. We will abuse the notation $\iota_n$ by also using it for the composition of the map above with the equality \eqref{eqn:k of stack}
\begin{equation}
\label{eqn:iota n}
K_{\BC^* \times \BC^*} (\comm_n) \xrightarrow{\iota_n} \BZ[q_1^{\pm 1}, q_2^{\pm 2}][z_1^{\pm 1}, \dots, z_n^{\pm 1}]^{\sym}
\end{equation}

\medskip

\begin{definition}
\label{big shuffle}

(\cite{FHHSY}) Consider the rational function
$$
\zeta(x) = \frac {(1-xq_1)(1-xq_2)\left(1-\frac {q_1q_2}{x} \right)}{1-x}
$$ 
The vector space
\begin{equation}
\label{eqn:big shuf}
\CV = \bigoplus_{n=0}^{\infty} \BZ[q_1^{\pm 1}, q_2^{\pm 2}][z_1^{\pm 1}, \dots, z_n^{\pm 1}]^{\esym}
\end{equation}
is made into an algebra via the following \textbf{shuffle product}
\begin{multline}
\label{eqn:shuf prod}
R (z_1,\dots,z_n) * R'(z_1,\dots,z_{n'}) = \\ = \eSym \left[ \frac {R (z_1,\dots,z_n) R'(z_{n+1},\dots,z_{n+n'})}{n! n'!} \prod_{1 \leq i \leq n < j \leq n+n'} \zeta \left(\frac {z_i}{z_j} \right) \right]
\end{multline}
(above, ``\emph{Sym}" refers to symmetrization with respect to $z_1,\dots,z_{n+n'}$).

\end{definition}

\medskip

\noindent The \textbf{$K$-theoretic Hall algebra} is defined as
$$
K = \bigoplus_{n=0}^{\infty} K_{\BC^* \times \BC^*} (\comm_n)
$$
It is endowed with a certain convolution product (\cite{SV Hilb}, see \cite[Section 2.3]{Quiver 1} for the construction in notation closer to ours), which has the property that the maps \eqref{eqn:iota n} combine to an algebra homomorphism
\begin{equation}
\label{eqn:iota}
K \xrightarrow{\iota} \CV
\end{equation}
Unfortunately, we do not know how to effectively describe the image of $\iota$.

\medskip

\subsection{}

In the present paper, we will study a variant of the $K$-theoretic Hall algebra, which we will be able to describe completely in terms of (the natural analogue of) the homomorphism \eqref{eqn:iota}. Consider the \textbf{semi-nilpotent commuting variety}
$$
\ocomm^{\nilp}_n \subset \ocomm_n
$$ 
parametrizing those pairs $(X,Y)$ of commuting $n \times n$ matrices with $X$ arbitrary and $Y$ nilpotent. The semi-nilpotency condition initially arose in the context of $K$-theoretic Hall algebras in \cite{SV gen}, but it also naturally arises in categorification via knot invariants (\cite{GN, GNR}). Letting the \textbf{semi-nilpotent commuting stack} be 
$$
\comm_n^{\nilp} = \ocomm_n^{\nilp}/GL_n
$$
we may define the following analogue of the construction of the previous Subsection
$$
K^{\nilp} = \bigoplus_{n=0}^{\infty} K_{\BC^* \times \BC^*} (\comm^{\nilp}_n)
$$
Then we have the natural analogue of the map \eqref{eqn:iota}
\begin{equation}
\label{eqn:iota nilp}
K^{\nilp} \xrightarrow{\iota^\nilp} \CV
\end{equation}
One endows $K^{\nilp}$ with the same kind of convolution product as $K$, thus making \eqref{eqn:iota nilp} into an algebra homomorphism. The map \eqref{eqn:iota nilp} is well-known to be injective (\cite[Lemma 2.5.1]{VV}). The main purpose of the present paper is to explicitly and effectively describe its image. We will actually provide two descriptions of the image: one as an explicit subalgebra $\CS \subset \CV$ (in Section \ref{sec:shuf}) and one by producing an explicit PBW basis of $\CS$ over the ring $\BZ[q_1^{\pm 1}, q_2^{\pm 1}] = K_{\BC^* \times \BC^*}(\text{pt})$ (in Section \ref{sec:pbw}).

\medskip

\subsection{} 

For any partition $\lambda = (n_1 \geq \dots \geq n_k) \vdash n$, let
\begin{equation}
\label{eqn:components}
\comm^{\nilp}_\lambda \stackrel{i_\lambda}\hookrightarrow \comm^{\nilp}_n
\end{equation}
be the closure of the substack consisting of pairs of commuting matrices $(X,Y)$ which (up to conjugation) are block triangular with respect to a flag of subspaces 
$$
0 = V_0 \subset V_1 \subset \dots \subset V_{k-1} \subset V_k = \BC^n
$$
where $\dim V_i/V_{i-1} = n_i$; above, ``block-triangular" means that
\begin{equation}
\label{eqn:block}
X(V_i) \subset V_i \qquad \text{and} \qquad Y(V_i) = V_{i-1}
\end{equation}
for all $i \in \{1,\dots,k\}$. A well-known fact of linear algebra is that
$$
\comm^{\nilp}_n = \bigcup_{\lambda \vdash n} \comm^{\nilp}_\lambda
$$
The substack corresponding to $\lambda = (n)$ is simply $\BA^{n^2}/GL_n$, as \eqref{eqn:block} requires $X$ to be arbitrary but $Y$ to be 0. As such, the composition
$$
K_{G_n}(\BA^{n^2}) \xrightarrow{i_{(n)*}} K_{G_n}(\ocomm^{\nilp}_n) \longrightarrow K_{G_n}(\BA^{2n^2}) \stackrel{|_\circ}\cong K_{G_n} (\text{pt})
$$
is simply given by mapping $X$ to $(X,0)$ and then restricting to the origin. Because of this, the image of the composition above is the principal ideal generated by the (equivariant) Koszul complex of
$$
\BA^{n^2} \hookrightarrow \BA^{2n^2}, \qquad X \mapsto (X,0)
$$
in the ring $\BZ[q_1^{\pm 1}, q_2^{\pm 1}][z_1^{\pm 1}, \dots, z_n^{\pm 1}]^{\sym}$. This Koszul complex is none other than
$$
F_n(z_1,\dots,z_n) = \prod_{1\leq i , j \leq n} \left(1 - \frac {z_iq_2}{z_j} \right) 
$$
Therefore, we have for all $n \in \BN$
\begin{equation}
\label{eqn:generate}
\BZ[q_1^{\pm 1}, q_2^{\pm 1}][z_1^{\pm 1}, \dots, z_n^{\pm 1}]^{\sym} \cdot F_n \subset \text{Im } \iota^{\nilp}
\end{equation}

\medskip

\begin{proposition}
\label{prop:generate}

As a $\BZ[q_1^{\pm 1}, q_2^{\pm 1}$]-algebra, $\emph{Im }\iota^{\enilp}$ is generated by the elements in the left-hand side of \eqref{eqn:generate}, as $n$ goes over $\BN$. 

\end{proposition}

\medskip 

\noindent The result above was proved at the level of Chow groups in \cite[Proposition 5.12]{SV gen}; the adaptation of the proof of \loccit to $K$-theory is straightforward, so we leave it as an exercise to the reader.

\bigskip

\section{The shuffle algebra}
\label{sec:shuf}

\medskip

\subsection{} The main purpose of the present Section is to identify the image of the map \eqref{eqn:iota nilp}. Proposition \ref{prop:generate} implies that
\begin{equation}
\label{eqn:image is generated}
\text{Im } \iota^{\nilp} = \Big \langle \BZ[q_1^{\pm 1}, q_2^{\pm 1}][z_1^{\pm 1}, \dots, z_n^{\pm 1}]^{\sym} \cdot F_n\Big \rangle_{n \in \BN} \subset \CV
\end{equation}

\medskip

\begin{definition}
\label{def:shuffle}

Consider the vector subspace $\CS \subset \CV$ consisting of symmetric Laurent polynomials $R(z_1,\dots,z_n)$ such that for any partition $(n_1 \geq \dots \geq n_k) \vdash n$, 
\begin{equation}
\label{eqn:specialization}
R(x_1,x_1q_2,\dots,x_1q_2^{n_1-1}, \dots, x_k,x_kq_2,\dots,x_k q_2^{n_k-1})
\end{equation}
is divisible by
\begin{multline}
\label{eqn:divisible by}
 \prod_{i=1}^k \left[ (1-q_2)^{n_i} \prod_{s=1}^{n_i-1} \zeta(q_2^s)^{n_i-s} \right] \\ \prod_{1 \leq i < j \leq k} \left[ \prod_{a = 1}^{n_i-1} \prod_{b=0}^{n_j-1} (x_i q_1 - x_j q_2^{b-a}) \right] \left[ \prod_{a = 1}^{n_i-1} \prod_{b=0}^{n_j-1} (x_j q_1 - x_i q_2^{a-b-1}) \right]
\end{multline}
in the ring $\BZ[q_1^{\pm 1}, q_2^{\pm 1}][x^{\pm 1}_1,\dots,x^{\pm 1}_k]$. We will call $\CS$ the (integral) \textbf{shuffle algebra}.

\end{definition}

\medskip

\begin{remark} Upon tensoring with $\BQ(q_1,q_2)$, all scalars $1-q_2$ and $\zeta(q_2^x)$ become invertible, and the fact that the specialization \eqref{eqn:specialization} is divisible by \eqref{eqn:divisible by} reduces to
\begin{equation}
\label{eqn:wheel}
R(x,xq_2,xq_1q_2,z_4,\dots,z_n) = R(x,xq_2,xq_1^{-1},z_4,\dots,z_n) = 0
\end{equation}
(which is none other than the particular case of \eqref{eqn:divisible by} for the partition $(2,1,\dots,1)\vdash n$). Conditions \eqref{eqn:wheel} are precisely the well-known wheel conditions (\cite{FHHSY}) for the shuffle algebra associated to quantum toroidal $\fgl_1$ over the field $\BQ(q_1,q_2)$.

\end{remark}

\medskip

\begin{remark} In the context of integral forms of quantum affine groups, divisibility conditions on integral shuffle algebras were first studied in \cite[Definition 3.37]{Ts}.

\end{remark}

\medskip

\subsection{}

The following two Propositions immediately establish the fact that 
\begin{equation}
\label{eqn:one-sided}
\text{Im } \iota^{\nilp} \subseteq \CS
\end{equation}

\medskip

\begin{proposition}
\label{prop:1}

For any $n \in \BN$, we have $\BZ[q_1^{\pm 1}, q_2^{\pm 1}][z_1^{\pm 1}, \dots, z_n^{\pm 1}]^{\esym} \cdot F_n \subset \CS$.

\end{proposition}

\medskip

\begin{proof} Since $F_n$ vanishes whenever we set $z_i = q_2z_j$ (for any $i \neq j$) then the conditions in Definition \ref{def:shuffle} for any multiple of $F_n$ are trivially satisfied.

\end{proof}

\begin{proposition}
\label{prop:2}

The subspace $\CS \subset \CV$ is a subalgebra with respect to \eqref{eqn:shuf prod}.

\end{proposition}

\medskip

\begin{proof} Let us write $\CS_n \subset \CS$ for the graded part consisting of Laurent polynomials in $n$ variables (i.e. the $n$-th direct summand of \eqref{eqn:big shuf}). We need to prove that 
$$
R \in \CS_{n'} \quad \text{and} \quad R' \in \CS_{n-n'} \qquad \Rightarrow \qquad R * R' \in \CS_{n}
$$
Consider the specialization of the set of variables $\{z_1,\dots,z_{n}\}$ at
\begin{equation}
\label{eqn:spec shuf}
\Big\{ x_1, x_1 q_2, \dots, x_1q_2^{n_1-1}, \dots, x_k, x_k q_2, \dots, x_kq_2^{n_k-1} \Big\}
\end{equation}
for some $n_1+\dots+n_k = n$. By \eqref{eqn:shuf prod}, to plug this specialization into $R * R'$ means to sum over all ways to permute the variables \eqref{eqn:spec shuf} and to plug them into
\begin{equation}
\label{eqn:second line}
 \frac {R (z_1,\dots,z_{n'}) R'(z_{n'+1},\dots,z_{n})}{n'! (n-n')!} \prod_{1 \leq i \leq n' < j \leq n} \zeta \left(\frac {z_i}{z_j} \right)
\end{equation}
However, because $\zeta(q_2^{-1}) = 0$, such a permutation can produce a non-zero contribution only if the variables
\begin{align*}
&x_i,x_iq_2,\dots,x_iq_2^{m_i-1} \qquad \quad \text{ are plugged into the variables of }R' \\
&x_iq_2^{m_i} ,x_iq_2^{m_i+1},\dots,x_iq_2^{n_i-1} \text{ are plugged into the variables of }R
\end{align*}
for some $m_i \in \{0,\dots,n_i\}$, for all $i \in \{1,\dots,k\}$. The contribution of such a permutation to the specialization of \eqref{eqn:second line} is then
$$
R(\dots,x_iq_2^{m_i} ,x_iq_2^{m_i+1},\dots,x_iq_2^{n_i-1}, \dots) R'(\dots,x_i,x_iq_2,\dots,x_iq_2^{m_i-1}, \dots) 
$$
\begin{equation}
\label{eqn:expression}
\prod_{i = 1}^k \prod_{a=m_i}^{n_i-1} \prod_{b=0}^{m_i-1} \zeta(q_2^{a-b}) \left[  \prod_{1 \leq i \neq j \leq k} \prod_{a=m_i}^{n_i-1} \prod_{b=0}^{m_j-1} \zeta \left(\frac {x_iq_2^a}{x_j q_2^b} \right)  \right]
\end{equation}
It remains to show that, for any $m_i \in \{0,\dots,n_i\}$, the expression \eqref{eqn:expression} is divisible by \eqref{eqn:divisible by}. Because $R \in \CS_{n'}$ and $R' \in \CS_{n-n'}$, the first line of \eqref{eqn:expression} is divisible by
$$
(1-q_2)^{n} \prod_{i=1}^k \left[ \prod_{s=1}^{m_i-1} \zeta(q_2^s)^{m_i-s} \prod_{s=1}^{n_i-m_i-1} \zeta(q_2^s)^{n_i-m_i-s}  \right]
$$
Together with the various $\zeta(q_2^{a-b})$ on the second line of \eqref{eqn:expression}, this precisely establishes divisibility by the expression on the first line of \eqref{eqn:divisible by}. Then it remains to prove that \eqref{eqn:expression} is divisible by the expression on the second line of \eqref{eqn:divisible by}. To this end, note that the formula in square brackets in \eqref{eqn:expression} is divisible by
\begin{equation}
\label{eqn:prod 1}
\prod_{1\leq i \neq j \leq k} \left[ \prod_{a=m_i}^{n_i-1} \prod_{b=0}^{m_j-1} (x_iq_1 - x_j q_2^{b-a})\prod_{a=0}^{m_i-1} \prod_{b=m_j}^{n_j-1} (x_i q_1 - x_j q_2^{b-a-1}) \right]
\end{equation}
Meanwhile, the second line of \eqref{eqn:divisible by} can be rewritten in a more symmetric way as
\begin{equation}
\label{eqn:prod final}
\prod_{1\leq i \neq j \leq k} \left[ \prod_{a = \min(n_i-n_j,0)+1}^{n_i-1} \prod_{b=0}^{\min(n_i,n_j)-1} (x_i q_1 - x_j q_2^{b-a}) \right]
\end{equation}
As a consequence, Definition \ref{def:shuffle} implies that the first line of \eqref{eqn:expression} is divisible by
\begin{multline}
\prod_{1\leq i \neq j \leq k} \left[  \prod_{a = \min(m_i-m_j,0)+1}^{m_i-1} \prod_{b=0}^{\min(m_i,m_j)-1} (x_i q_1 - x_j q_2^{b-a}) \right. \\
\left.  \prod_{a = \min(m_i-m_j,n_i-n_j)+1}^{n_i-m_j-1} \prod_{b=0}^{\min(n_i-m_i,n_j-m_j)-1} (x_i q_1 - x_j q_2^{b-a}) \right] \label{eqn:prod 2}
\end{multline}
It is elementary \footnote{Explicitly, this ``elementary" claim follows from the fact that the Laurent polynomial
$$
\sum_{a=m_i}^{n_i-1} \sum_{b=0}^{m_j-1} z^{b-a} + \sum_{a=0}^{m_i-1} \sum_{b=m_j}^{n_j-1} z^{b-a-1} + \sum_{a = \min(m_i-m_j,0)+1}^{m_i-1} \sum_{b=0}^{\min(m_i,m_j)-1} z^{b-a} +
$$
$$
+ \sum_{a = \min(m_i-m_j,n_i-n_j)+1}^{n_i-m_j-1} \sum_{b=0}^{\min(n_i-m_i,n_j-m_j)-1} z^{b-a} -  \sum_{a = \min(n_i-n_j,0)+1}^{n_i-1} \sum_{b=0}^{\min(n_i,n_j)-1} z^{b-a}
$$
has non-negative coefficients, as it is equal to $\sum_{a=\min(n_j-n_i,m_j-m_i)}^{\max(0,m_j-m_i)-1} z^a$.} to see that the product of \eqref{eqn:prod 1} and \eqref{eqn:prod 2} is divisible by \eqref{eqn:prod final}, for any choice of numbers $m_i \in \{0,\dots,n_i\}$, exactly what we needed to prove.

\end{proof}

\medskip

\subsection{} We will now prove the opposite inclusion to \eqref{eqn:one-sided}, thus concluding the proof of Theorem \ref{thm:main}.

\medskip

\begin{proposition}
\label{prop:magic}

We have $\emph{Im } \iota^{\enilp} \supseteq \CS$.

\end{proposition}

\medskip

\begin{proof} We will refine the argument of \cite[Proposition 2.4]{Shuf}, itself based on \cite{FHHSY}. For any partition $\lambda = (n_1 \geq \dots \geq n_k) \vdash n$, consider the linear map
$$
\CS_n \xrightarrow{\ph_{\lambda}} \BZ[q_1^{\pm 1}, q_2^{\pm 1}][x_1^{\pm 1}, \dots ,x_k^{\pm 1}]
$$
\begin{equation}
\label{eqn:spec map}
R(z_1,\dots,z_n) \mapsto R(x_1q_2^{n_1-1}, \dots, x_1q_2,x_1, \dots, x_k q_2^{n_k-1}, \dots, x_kq_2, x_k)
\end{equation}
We consider the total lexicographic order on partitions of size $n$, where 
$$
(m_1 \geq \dots \geq m_l) > (n_1 \geq \dots \geq  n_k)
$$
means that there exists $i$ such that $m_1 = n_1, \dots, m_{i} = n_{i}, m_{i+1} > n_{i+1}$. The sets
$$
\CS_\lambda = \bigcap_{\mu > \lambda} \text{Ker } \ph_{\mu}
$$
yield an increasing filtration of $\CS_n = \CS_{(n)}$.

\medskip

\begin{claim}
\label{claim:magic}

For any $R \in \CS_{\lambda}$, there exists $R' \in (\emph{Im } \iota^{\enilp})\cap \CS_{\lambda}$ such that
\begin{equation}
\label{eqn:phis equal}
\ph_{\lambda}(R) = \ph_{\lambda}(R')
\end{equation}

\end{claim}

\medskip

\noindent Iterating Claim \ref{claim:magic} for all partitions $\lambda$ in decreasing lexicographic order allows us to take any $R \in \CS_n$, and by subtracting various elements in $\text{Im } \iota^{\nilp}$, ensure that it lies in the kernel of $\ph_{\lambda}$ for smaller and smaller $\lambda$. As soon as we pass $\lambda = (1,\dots,1)$, then we will have made $R$ equal to 0 by subtracting various elements in $\text{Im } \iota^{\nilp}$, and the proof of Proposition \ref{prop:magic} would be complete.

\medskip

\noindent Let us now prove Claim \ref{claim:magic}. If we write $\lambda = (n_1 \geq \dots \geq n_k)$, then the transposed partition $\lambda' = (t_1 \geq \dots \geq t_p)$ is defined by the equation
\begin{equation}
\label{eqn:connection transpose}
n_i = |u\in \{1,\dots,p\} \text{ s.t. } t_u \geq i| 
\end{equation}
for all $i$. Let us write $s_i = t_1+\dots+t_i$ for all $i$, and define
\begin{equation}
\label{eqn:that r'}
R'(z_1,\dots,z_n) = \Sym \Big[ r(z_1,\dots,z_n) \Big]
\end{equation}
where for any $\rho \in \BZ[q_1^{\pm 1}, q_2^{\pm 1}][z_1^{\pm 1}, \dots, z_n^{\pm 1}]^{\text{psym}}$, we set
$$
r = \rho(z_1,\dots,z_n) \prod_{i=1}^p F_{t_i}(z_{s_{i-1}+1},\dots,z_{s_i}) \prod_{1 \leq i < j \leq p} \prod_{a=s_{i-1}+1}^{s_i} \prod_{b = s_{j-1}+1}^{s_j} \zeta \left(\frac {z_a}{z_b} \right)
$$
(the superscript ``psym" means that we require $\rho$ to be symmetric in $z_{s_{i-1}+1},\dots,z_{s_i}$ for all $i \in \{1,\dots, p\}$ separately). We claim that
\begin{equation}
\label{eqn:as it should}
R' \in (\text{Im } \iota^{\nilp}) \cap \CS_{\lambda} 
\end{equation}
It is easy to see that $R' \in \text{Im } \iota^{\nilp}$, as it equals the shuffle product of Laurent polynomials divisible by $F_{t_1}, \dots, F_{t_p}$. Moreover, we claim that $R' \in \CS_{\lambda}$; to see this, we must show that $R'$ is annihilated by $\ph_{\mu}$ for any $\mu > \lambda$. Indeed, computing $\ph_{\mu}(R')$ for some $\mu = (m_1 \geq \dots \geq m_l)$ entails specializing the variables of $R'$ to 
\begin{equation}
\label{eqn:we specialize}
\Big\{ x_i, x_i q_2, \dots ,x_iq_2^{m_i-1} \Big\}_{\{1,\dots,l\}}
\end{equation}
Equivalently, this amounts to inserting the variables \eqref{eqn:we specialize} among the arguments of $r$ in an arbitrary order. Let's call such an insertion ``good" if for each $i \in \{1,\dots,l\}$, the variables $x_iq_2^{m_i-1},\dots,x_iq_2,x_i$ are plugged in successive sets among
\begin{equation}
\label{eqn:chunks}
\Big\{z_1,\dots,z_{s_1}\Big\}, \Big\{z_{s_1+1},\dots,z_{s_2}\Big\}, \dots, \Big\{z_{s_{p-1}+1}, \dots, z_n\Big\}
\end{equation}
Because $\zeta(q_2^{-1}) = 0$ and $F_t(\dots,x,xq_2,\dots) = 0$ for all $t$, only good insertions have the property that $r$ specializes to a non-zero value. However, $\mu > \lambda$ means that
\begin{align*}
&m_1 = n_1 = |u\in \{1,\dots,p\} \text{ s.t. } t_u \geq 1| \\ 
&\dots \\
&m_i = n_i = |u\in \{1,\dots,p\} \text{ s.t. } t_u \geq i|  \\
&m_{i+1} > n_{i+1} = |u\in \{1,\dots,p\} \text{ s.t. } t_u \geq i+1|  
\end{align*}
for some $i$, and thus good insertions cannot exist. This establishes \eqref{eqn:as it should}. 

\medskip

\noindent It remains to show that we can choose the Laurent polynomial $\rho$ in the definition of $r$ so that \eqref{eqn:phis equal} holds. Repeating the argument in the preceding paragraph shows that $\ph_{\lambda}(R')$ is calculated by inserting the variables $x_iq_2^{n_i-1}, \dots, x_iq_2,x_i$ in the arguments of $r$. In this case, the only good insertions are those such that
$$
\Big\{ z_{s_{i-1}+1}, z_{s_{i-1}+2}, \dots, z_{s_i} \Big\} = \Big\{ x_1q_2^{n_1-i}, x_2q_2^{n_2-i}, \dots, x_{t_i} q_2^{n_{t_i}-i} \Big\}
$$
for all $i \in \{1,\dots,p\}$. Thus, we conclude that
$$
\ph_{\lambda}(R') =  \ph_{\lambda} (\rho)  \prod_{i=1}^p \prod_{1\leq a, b \leq t_i} \left(1- \frac {x_aq_2^{n_a+1}}{x_bq_2^{n_b}}\right) \prod_{1\leq i < j \leq p} \prod_{a=1}^{t_i} \prod_{b = 1}^{t_j} \zeta \left(\frac {x_aq_2^{n_a-i}}{x_bq_2^{n_b-j}} \right) 
$$
Although $\rho$ is not itself an element of $\CS_n$, the notation $\ph_{\lambda} (\rho)$ is defined just like \eqref{eqn:spec map}. We may now move the products in $a,b$ from the inside to the outside of the above formula, and obtain (after clearing various cancelations involving $\zeta$ factors)
\begin{equation}
\label{eqn:final}
\ph_{\lambda}(R') = \ph_{\lambda} (\rho) \cdot \Pi_1 \Pi_2 \Pi_3
\end{equation}
where
\begin{align*}
&\Pi_1 = \prod_{a=1}^k \left[ (1-q_2)^{n_a} \prod_{u=1}^{n_a-1} \zeta(q_2^u)^{n_a-u} \right] \\
&\Pi_2 = \prod_{1\leq a \neq b \leq k} \left[ \mathop{\prod_{0 \leq u < n_a}^{u-v > n_a-n_b}}_{0 \leq v < n_b} \left(1 - \frac {x_aq_1}{x_b q_2^{v-u}} \right) \mathop{\prod_{1 \leq u \leq n_a}^{u-v \leq n_a-n_b}}_{0 \leq v < n_b} \left(1 - \frac {x_aq_1}{x_b q_2^{v-u}} \right) \right] \\
&\Pi_3 = \prod_{1\leq a \neq b \leq k} \prod_{u=\max(n_a-n_b,0)+1}^{n_a} \left(1 - \frac {x_aq_2^u}{x_b} \right)
\end{align*}
Clearly, $\Pi_1$ is precisely the first line of \eqref{eqn:divisible by}, while it is elementary to see that $\Pi_2$ matches the second line of \eqref{eqn:divisible by} up to an overall monomial. The fact that $R \in \CS_n$ implies that $\ph_{\lambda}(R)$ is divisible by \eqref{eqn:divisible by}, and thus is divisible by $\Pi_1\Pi_2$. However, the fact that $R \in \CS_{\lambda}$ implies certain additional divisibilities: whenever
$$
x_i q_2^{-1} \quad \text{or} \quad x_iq_2^{n_i} \qquad \text{is set equal to} \qquad x_j,x_j q_2,\dots, x_j q_2^{n_j-1}
$$
for some $i<j$, the quantity $\ph_{\lambda}(R)$ must vanish (indeed, this is because if we enlarge $n_i$ and diminish $n_j$ by some positive amount, the resulting partition $\mu$ is larger than $\lambda$). This precisely entails the fact that $\ph_{\lambda}(R)$ is divisible by $\Pi_3$, so we conclude that there exists a Laurent polynomial $A(x_1,\dots,x_k)$ such that
\begin{equation}
\label{eqn:final 2}
\ph_{\lambda}(R) = A(x_1,\dots,x_k) \cdot \Pi_1\Pi_2\Pi_3
\end{equation}
Moreover, $A(x_1,\dots,x_k)$ is symmetric in $x_a$ and $x_b$ if $n_a=n_b$, because $R$ is symmetric in all of its variables. Thus, we must choose $\rho$ such that
\begin{equation}
\label{eqn:psym =}
\ph_{\lambda} (\rho) = A(x_1,\dots,x_k)
\end{equation}
and then \eqref{eqn:final} and \eqref{eqn:final 2}  would imply \eqref{eqn:phis equal}. We may assume that $A$ is a polynomial in $x_1,\dots,x_k$, by multiplying \eqref{eqn:psym =} with a sufficiently high monomial. Thus, if the partition $\lambda$ consists of $d_1$ times $1$, $d_2$ times $2$ etc, we may assume that
$$
A(x_1,\dots,x_k) = m_{\nu_1}(x_k, \dots,x_{k-d_1+1}) m_{\nu_2}(x_{k-d_1}, \dots, x_{k-d_1-d_2+1}) \dots
$$
where $m_{\nu}(z_1,z_2,\dots) = \Sym [z_1^{\nu_1} z_2^{\nu_2} \dots]$ denotes the monomial symmetric function associated to the partition $\nu = (\nu_1 \geq \nu_2 \geq \dots)$. If we define
$$
\rho'(z_1,\dots,z_n) = m_{\nu_1}(z_1,\dots,z_{t_1}) m_{\nu_2}(z_{t_1+1}, \dots, z_{t_1+t_2}) \dots
$$
then it is straightforward to see that
$$
\ph_{\lambda} (\rho') = q_2^{\text{some integer}} \cdot A(x_1,\dots,x_k) + B(x_1,\dots,x_k)
$$
where $B$ is a polynomial, symmetric in $x_a$ and $x_b$ if $n_a=n_b$, for which the sequence
$$
\left( \underset{x_k, \dots, x_{k-d_1}}{\text{hom deg}} B, \underset{x_{k-d_1+1}, \dots, x_{k-d_1-d_2}}{\text{hom deg}} B, \dots \right)
$$
is lexicographically smaller than the analogous sequence for $A$. Therefore, we may repeat the argument above for $B$ instead of $A$; after finitely many iterations of this procedure, we would obtain a polynomial $\rho$ for which \eqref{eqn:psym =} holds precisely.

\end{proof}

\medskip

\begin{proof} \emph{of Theorem \ref{thm:main}:} immediate from \eqref{eqn:one-sided} and Proposition \ref{prop:magic}.

\end{proof}

\bigskip

\section{The PBW basis}
\label{sec:pbw}

\medskip

\subsection{} For any $(n,d) \in \BN \times \BZ$, consider the Laurent polynomial
\begin{equation}
\label{eqn:pnd}
P_{n,d} = \Sym \left[\frac {\prod_{i=1}^n z_i^{\left \lfloor \frac {id}n \right \rfloor - \left \lfloor \frac {(i-1)d}n \right \rfloor} \sum_{s=0}^{t-1} \frac {z_{a(t-1)+1} \dots z_{a(t-s)+1}}{q_2^s z_{a(t-1)} \dots z_{a(t-s)}}}{\prod_{i=1}^{n-1} \left(1-\frac {z_{i+1}}{z_iq_2}\right)} \prod_{1\leq i < j \leq n} \zeta \left(\frac {z_i}{z_j} \right) \right]
\end{equation}
where we write $t = \gcd(n,d)$ and $a = \frac nt$. With the notation above, let
\begin{equation}
\label{eqn:gamma}
\gamma_{n,d} = \frac {q_2^t - 1}{(q_1^t - 1)(q_3^t - 1)} \cdot (q_1-1)^n (q_3 - 1)^n
\end{equation}
with $q_3 = \frac 1{q_1q_2}$. Let us define the following rescaled versions of \eqref{eqn:pnd}
\begin{equation}
\label{eqn:pnd rescaled}
\oP_{n,d} = \gamma_{n,d} \cdot P_{n,d}
\end{equation}
A sequence $v = \{(n_1,d_1), \dots, (n_k,d_k)\} \subset \BN \times \BZ$ will be called a \textbf{convex path} if
$$
\frac {d_1}{n_1} \leq \dots \leq \frac {d_k}{n_k}
$$
We always consider convex paths up to the equivalence generated by permuting lattice points of the same slope. This is motivated by the fact that $P_{n,d} $ and $P_{n',d'}$ commute if $(n,d)$ and $(n',d')$ have the same slope (\cite{Shuf}), and thus the expressions
\begin{align}
&P_v = P_{n_1,d_1} * \dots * P_{n_k,d_k} \label{eqn:path 1} \\
&\oP_v = \oP_{n_1,d_1} * \dots * \oP_{n_k,d_k}  \label{eqn:path 2}
\end{align}
only depend on the equivalence class of a convex path. It was shown in \cite{Shuf} that
\begin{equation}
\label{eqn:pbw loc}
\CS_{\loc} := \CS \bigotimes_{\BZ[q_1^{\pm 1}, q_2^{\pm 1}]} \BQ(q_1,q_2) = \bigoplus_{v \text{ convex path}} \BQ(q_1,q_2) \cdot P_v
\end{equation}
following the analogous result of \cite{BS} for the elliptic Hall algebra. Because relation \eqref{eqn:pbw loc} is taken over $\BQ(q_1,q_2)$, it also holds with the $P$'s replaced by $\oP$'s. 

\medskip

\begin{remark}

The following formulas are proved in \cite[(2.34) and (2.35)]{W-alg}
\begin{align*}
P_{n,d} &= \gamma'_{n,d} \cdot \eSym \left[\frac {\prod_{i=1}^n z_i^{\left \lfloor \frac {id}n \right \rfloor - \left \lfloor \frac {(i-1)d}n \right \rfloor} \sum_{s=0}^{t-1} \frac {z_{a(t-1)+1} \dots z_{a(t-s)+1}}{q_1^s z_{a(t-1)} \dots z_{a(t-s)}}}{\prod_{i=1}^{n-1} \left(1-\frac {z_{i+1}}{z_iq_1}\right)} \prod_{1\leq i < j \leq n} \zeta \left(\frac {z_i}{z_j} \right) \right] \\
&= \gamma''_{n,d} \cdot \eSym \left[\frac {\prod_{i=1}^n z_i^{\left \lfloor \frac {id}n \right \rfloor - \left \lfloor \frac {(i-1)d}n \right \rfloor} \sum_{s=0}^{t-1} \frac {  z_{a(t-1)+1} \dots z_{a(t-s)+1}}{q_3^s z_{a(t-1)} \dots z_{a(t-s)}}}{\prod_{i=1}^{n-1} \left(1-\frac {z_{i+1}}{z_iq_3}\right)} \prod_{1\leq i < j \leq n} \zeta \left(\frac {z_i}{z_j} \right) \right] 
\end{align*}
where we recall that $q_3 = \frac 1{q_1q_2}$, $t = \gcd(n,d)$, $a = \frac nt$, and define
$$
\gamma'_{n,d} = \frac {q_1^t-1}{(q_1-1)^n} \cdot \frac {(q_2-1)^n}{q_2^t-1} \qquad \text{and} \qquad \gamma''_{n,d} = \frac {q_3^t-1}{(q_3-1)^n} \cdot \frac {(q_2-1)^n}{q_2^t-1}
$$

\end{remark}

\medskip

\subsection{} For any coprime $(n,d) \in \BN \times \BZ$, the following power series identities
\begin{align}
&1 + \sum_{t=1}^{\infty} \frac {H_{nt,dt}}{x^t} = \exp \left(\sum_{t=1}^{\infty} \frac {P_{nt,dt}}{tx^t} \right) \label{eqn:h to p} \\ 
&1 + \sum_{t=1}^{\infty} \frac {\oH_{nt,dt}}{x^t} = \exp \left(\sum_{t=1}^{\infty} \frac {\oP_{nst,dt}}{tx^t} \right) \label{eqn:oh to op}
\end{align}
define elements $\{H_{n,d}, \oH_{n,d}\}_{(n,d) \in \BN \times \BZ} \in \CS_{\loc}$. In \cite[Formula (2.9)]{AGT}, we showed that
\begin{equation}
\label{eqn:end}
H_{n,d} = \Sym \left[\frac {\prod_{i=1}^n z_i^{\left \lfloor \frac {id}n \right \rfloor - \left \lfloor \frac {(i-1)d}n \right \rfloor}}{\prod_{i=1}^{n-1} \left(1-\frac {z_{i+1}}{z_iq_2}\right)} \prod_{1\leq i < j \leq n} \zeta \left(\frac {z_i}{z_j} \right) \right]
\end{equation}
$\forall (n,d) \in \BN \times \BZ$. Similarly, the following formula can be found in \cite[Exercise 3.18]{Pieri}
\begin{multline}
\label{eqn:end 1}
\oH_{n,d} = (q_1-1)^n(q_2-1)^n \cdot \\ \Sym \left[\frac {\prod_{i=1}^n z_i^{\left \lfloor \frac {id}n \right \rfloor - \left \lfloor \frac {(i-1)d}n \right \rfloor} \prod_{s=1}^{t-1} \left(q_1^s - \frac {z_{as+1}}{z_{as} q_3} \right)}{(q_1-1)\dots(q_1^t-1) \prod_{i=1}^{n-1} \left(1-\frac {z_{i+1}}{z_iq_3}\right)} \prod_{1\leq i < j \leq n} \zeta \left(\frac {z_i}{z_j} \right) \right]  
\end{multline}
By switching the roles of $q_1$ and $q_3$, one also obtains the following analogous formula
\begin{multline}
\label{eqn:end 2}
\oH_{n,d} = (q_2-1)^n(q_3-1)^n \cdot \\ \Sym \left[\frac {\prod_{i=1}^n z_i^{\left \lfloor \frac {id}n \right \rfloor - \left \lfloor \frac {(i-1)d}n \right \rfloor} \prod_{s=1}^{t-1} \left(q_3^s - \frac {z_{as+1}}{z_{as} q_1} \right)}{(q_3-1)\dots(q_3^t-1) \prod_{i=1}^{n-1} \left(1-\frac {z_{i+1}}{z_iq_1}\right)} \prod_{1\leq i < j \leq n} \zeta \left(\frac {z_i}{z_j} \right) \right]  
\end{multline}

\medskip

\begin{proposition}
\label{prop:belong}

We have $\oH_{n,d} \in \CS$ for all $(n,d) \in \BN \times \BZ$.

\end{proposition}

\medskip

\begin{proof} Consider any partition $(n_1,\dots,n_k) \vdash n$. In the Laurent polynomial
$$
\oH_{n,d}(z_1,\dots,z_n) \cdot (1+q_3)(1+q_3+q_3^2)\dots(1+q_3+\dots+q_3^{t-1})
$$
(where $t = \gcd(n,d)$), let us specialize the variables $z_1,\dots,z_n$ to
\begin{equation}
\label{eqn:inserting}
x_1,x_1q_2,\dots,x_1q_2^{n_1-1}, \dots, x_k, x_kq_2, \dots, x_kq_2^{n_k-1}
\end{equation}
Using formula \eqref{eqn:end 2}, this amounts to permuting the variables \eqref{eqn:inserting} arbitrarily, and then inserting them instead of $z_1,\dots,z_n$ into a certain expression of the form
\begin{equation}
\label{eqn:zetas}
\frac {(1-q_2)^n \cdot \text{Laurent polynomial}}{\prod_{i=1}^{n-1} \left(1-\frac {z_{i+1}}{z_iq_1}\right)} \prod_{1\leq i < j \leq n} \zeta \left(\frac {z_i}{z_j} \right)
\end{equation}
Because $\zeta(q_2^{-1}) = 0$, the only insertions which produce a non-zero contribution are those for which $x_iq_2^{n_i-1}, \dots, x_i$ are plugged into $z_{a_1}, \dots, z_{a_{n_i}}$ for certain indices $a_1 < \dots < a_{n_i}$, for each $i \in \{1,\dots,k\}$. As such, it is clear that the resulting specialization is divisible by the expression on the first line of \eqref{eqn:divisible by}. Moreover, for any $i \neq j$ such that $n_i \geq n_j$, let us zoom in on a fixed $b \in \{0,\dots, n_j-1\}$ and assume that the variables \eqref{eqn:inserting} are permuted in the order
$$
x_i q_2^{n_i-1}, \dots, x_iq_2^u, x_jq_2^b, x_iq_2^{u-1}, \dots, x_i
$$
for some $u$. Then the product of $\zeta$ functions in \eqref{eqn:zetas} is a multiple of
\begin{multline*}
(x_iq_1-x_jq_2^{b-n_i+1}) \dots (x_i q_1-x_jq_2^{b-u-1}) \underline{(x_i q_1-x_jq_2^{b-u})} \\ (x_i q_1-x_jq_2^{b-u}) (x_i q_1-x_jq_2^{b-u+1}) \dots (x_i q_1-x_jq_2^{b-1})
\end{multline*}
As for the denominator in \eqref{eqn:zetas}, it can at most cancel the underlined term above. The resulting expression is a multiple of $\prod_{a=1}^{n_i-1} (x_iq_1 - x_j q_2^{b-a})$; taking the product over $b \in \{0,\dots,n_j-1\}$'s shows that the overall specialization is a multiple of the first product on the second line of \eqref{eqn:divisible by}. One shows that the specialization is a multiple of the second product on the second line of \eqref{eqn:divisible by} analogously. Thus
\begin{multline*}
\frac {\oH_{n,d}(x_1,\dots,x_1q_2^{n-1},\dots,x_k,\dots,x_kq_2^{n_k-1})}{\text{expression \eqref{eqn:divisible by}}} \in \\
\BZ[q_1^{\pm 1},q_2^{\pm 1}]_{(1+q_3+\dots+q_3^{s-1})_{s \in \BN}}[x_1^{\pm 1}, \dots,x_k^{\pm 1}]
\end{multline*}
Repeating the argument with the roles of $q_1$ and $q_3$ switched (i.e. using \eqref{eqn:end 1} instead of \eqref{eqn:end 2}) shows that the ratio above has coefficients in the localization
$$
\BZ[q_1^{\pm 1},q_2^{\pm 1}]_{(1+q_1+\dots+q_1^{s-1})_{s \in \BN}}
$$
We conclude that the coefficients are actually in $\BZ[q_1^{\pm 1},q_2^{\pm 1}]$, as we needed to show. 

\end{proof}

\medskip

\subsection{}

By analogy with \eqref{eqn:path 1}--\eqref{eqn:path 2}, let us write for any convex path $v$
\begin{align}
&H_v = H_{n_1,d_1} * \dots * H_{n_k,d_k} \label{eqn:path 3} \\
&\oH_v = \oH_{n_1,d_1} * \dots * \oH_{n_k,d_k}  \label{eqn:path 4}
\end{align}
Clearly, the elements $P_{n,d}$ may be replaced by either $H_{n,d}$ and $\oH_{n,d}$ in \eqref{eqn:pbw loc}, to produce a valid basis of $\CS_{\loc}$ as a $\BQ(q_1,q_2)$-vector space. However, our main interest is in the following $\BZ[q_1^{\pm 1}, q_2^{\pm 1}]$-submodule of $\CS_{\loc}$
$$
\CA = \bigoplus_{v \text{ convex path}} \BZ[q_1^{\pm 1}, q_2^{\pm 1}] \cdot \oH_v
$$
We are now ready to prove Lemma \ref{lem:pbw}, which provides an integral version of \eqref{eqn:pbw loc}. 

\medskip

\begin{proof} \emph{of Lemma \ref{lem:pbw}:} By Propositions \ref{prop:2} and \ref{prop:belong}, we have
$$
\CS \supseteq \CA
$$
It remains to prove the opposite inclusion, namely
\begin{equation}
\label{eqn:thing 1}
\CS \subseteq \CA
\end{equation}
To this end, recall the symmetric pairing defined in \cite[Formula (4.7)]{Shuf}
\begin{equation}
\label{eqn:the pairing}
\CS_{\loc} \otimes \CS_{\loc} \xrightarrow{\langle \cdot, \cdot \rangle} \BQ(q_1,q_2)
\end{equation}
by the formula \footnote{Note that our normalization of \eqref{eqn:formula pairing} differs from that of \loccit by $(q_1-1)^n (q_3-1)^n$.}
\begin{equation}
\label{eqn:formula pairing}
\Big \langle R, R' \Big \rangle = \frac 1{(q_2-1)^n} \int_{|z_1| \gg \dots \gg |z_n|} \frac {r(z_1,\dots,z_n)R (z_1^{-1}, \dots, z_n^{-1})}{\prod_{1\leq i < j \leq n} \zeta \left(\frac {z_j}{z_i} \right)} \prod_{i=1}^n Dz_i
\end{equation}
(where $Dz = \frac {dz}{2\pi i z}$), for any $R' \in \CS_{\loc}$ and
\begin{equation}
\label{eqn:rho}
R =  \Sym \left[r(z_1,\dots,z_n) \prod_{1 \leq i < j \leq n} \zeta \left(\frac {z_i}{z_j} \right) \right]
\end{equation}
where $r$ is an arbitrary Laurent polynomial with coefficients in $\BQ(q_1,q_2)$. Since any element $R \in \CS_{\loc}$ can be written in the form \eqref{eqn:rho} for some Laurent polynomial $r$ with coefficients in $\BQ(q_1,q_2)$ (as proved in \cite[Theorem 2.5]{Shuf}), formula \eqref{eqn:formula pairing} determines the pairing \eqref{eqn:the pairing} completely. It was shown in \cite[Proposition 5.7]{Shuf} that $\{P_v\}_{v \text{ convex}}$ is an orthogonal basis with respect to the pairing \eqref{eqn:the pairing}, satisfying 
\begin{equation}
\label{eqn:pair pv}
\Big \langle P_v, \oP_v \Big \rangle = \prod_{\mu \in \BQ} z_{\lambda_v^{\mu}}
\end{equation}
Let us explain the notation in the right-hand side of \eqref{eqn:pair pv}: for any convex path $v = \{(n_1,d_1), \dots, (n_k,d_k)\}$ and any 
$$
\mu = \frac dn \in \BQ
$$
(assume $\gcd(n,d) = 1$), those elements of $v$ of slope $\mu$ will be of the form
$$
(nt_1,dt_1), \dots,(nt_k,dt_k)
$$
for some partition $\lambda_v^{\mu} = (t_1 \geq \dots \geq t_k)$. As $\mu$ goes over the infinitely many rational numbers, all but finitely many of these partitions will be empty. Finally, for any partition $\lambda = (t_1 \geq \dots \geq t_k)$, we set
$$
z_\lambda = t_1\dots t_k \prod_{u \in \BN} (\text{number of }u\text{'s in }\lambda)!
$$
and this completes the explanation of the right-hand side of \eqref{eqn:pair pv}.

\medskip

\begin{claim}
\label{claim:sym}

An element $R \in \CS_{\eloc}$ is a linear combination of $\oH_v$'s with coefficients in $\BZ[q_1^{\pm 1}, q_2^{\pm 1}]$ if and only if
$$
\Big \langle R, H_v \Big \rangle \in \BZ[q_1^{\pm 1}, q_2^{\pm 1}]
$$
for all convex paths $v$. 

\end{claim}

\medskip 

\noindent Formula \eqref{eqn:pair pv} reduces Claim \ref{claim:sym} to the following well-known fact about symmetric functions: the Hall inner product of a symmetric function $f$ with all products of complete symmetric functions $h_n$ are integral if and only if $f$ is an integral linear combination of products of complete symmetric functions (indeed, products of complete symmetric functions yield the dual basis to monomial symmetric functions). 

\medskip 

\noindent Thus, Claim \ref{claim:sym} reduces \eqref{eqn:thing 1} to showing that 
\begin{equation}
\label{eqn:thing 2}
\Big \langle R, H_v \Big \rangle \in \BZ[q_1^{\pm 1}, q_2^{\pm 1}] 
\end{equation}
for any $R \in \CS$ and any convex path $v$. The remainder of the proof will deal with establishing \eqref{eqn:thing 2}. To this end, formula \eqref{eqn:end} implies that $H_v$ is a particular element of the shuffle algebra of the form
\begin{equation}
\label{eqn:r prime}
R'(z_1,\dots,z_n) = \Sym \left[\frac {p(z_1,\dots,z_n)}{\prod_{i=1}^{n-1} \left(1 - \frac {z_{i+1}}{z_iq_2}\right)} \prod_{1\leq i < j \leq n} \zeta \left( \frac {z_i}{z_j} \right) \right]
\end{equation}
where $p(z_1,\dots,z_n)$ is an arbitrary Laurent polynomial with coefficients in $\BZ[q_1^{\pm 1}, q_2^{\pm 1}]$.

\medskip

\begin{claim}
\label{claim:pairing}

For any $R \in \CS_{\emph{loc}}$ and any $R'$ as in \eqref{eqn:r prime}, we have
$$
\Big \langle R, R' \Big \rangle =  \sum_{n_1+\dots+n_k = n} \int_{|x_1| \ll \dots \ll |x_k|} \left[ \underset{\{z_{n_1+\dots+n_{i-1}+1},\dots,z_{n_1+\dots+n_{i}}\}_{1\leq i \leq k} = \{x_iq_2^{n_i-1},\dots,x_i\}_{1\leq i \leq k}}{\emph{Res}} \right. 
$$
\begin{equation}
\label{eqn:claim pairing}
\left. \frac 1{(q_2-1)^n} \cdot \frac {R(z_1,\dots,z_n) \cdot p(z_1^{-1}, \dots, z_n^{-1})}{\prod_{i=1}^{n-1} \left(1 - \frac {z_i}{z_{i+1}q_2} \right) \prod_{1\leq i < j \leq n} \zeta \left(\frac {z_i}{z_j}\right)} \right] \prod_{i=1}^k Dx_i
\end{equation}
where 
$$
\underset{\{z_1,\dots,z_n\} = \{xq_2^{n-1},\dots,x\}}{\emph{Res}}
$$
denotes the iterated residue first at $z_{n-1} = z_nq_2$, then at $z_{n-2} = z_nq_2^2$, $\dots$, finally at $z_1 = z_nq_2^{n-1}$, followed by relabeling the variable $z_n$ by $x$. 

\medskip

\end{claim}

\noindent Let us first indicate how Claim \ref{claim:pairing} implies \eqref{eqn:thing 2}. When $R \in \CS$, Definition \ref{def:shuffle} tells us that the expression in square brackets of \eqref{eqn:claim pairing} is a Laurent polynomial with coefficients in $\BZ[q_1^{\pm 1}, q_2^{\pm 1}]$, divided by various linear terms of the form
$$
x_i q_1^{\dots} - x_j q_2^{\dots}
$$
for $i \neq j$. As we take the integral of such an expression in the limit $|x_1| \ll \dots \ll |x_k|$, the result is still an element of $\BZ[q_1^{\pm 1}, q_2^{\pm 1}]$, thus establishing \eqref{eqn:thing 2}.

\medskip

\begin{proof} \emph{of Claim \ref{claim:pairing}:} It suffices to prove \eqref{eqn:claim pairing} for $R$ of the form \eqref{eqn:rho}, as such elements span $\CS_{\loc}$. Then the right-hand side of \eqref{eqn:formula pairing} is equal to $(q_2-1)^{-n}$ times
$$
\int_{|z_1| \gg \dots \gg |z_n|} \frac {r(z_1,\dots,z_n)}{\prod_{1\leq i < j \leq n} \zeta \left(\frac {z_j}{z_i} \right)} \sum_{\sigma \in S(n)} \frac {p(z^{-1}_{\sigma(1)}, \dots, z^{-1}_{\sigma(n)}) \prod_{1\leq i < j \leq n} \zeta \left(\frac {z_{\sigma(j)}}{z_{\sigma(i)}} \right)}{\prod_{i=1}^{n-1} \left(1 - \frac {z_{\sigma(i)}}{z_{\sigma(i+1)}q_2} \right)} \prod_{i=1}^n Dz_i 
$$
$$
= \int_{|z_1| \gg \dots \gg |z_n|} \sum_{\sigma \in S(n)} \frac {r(z_1,\dots,z_n) \cdot p(z^{-1}_{\sigma(1)}, \dots, z^{-1}_{\sigma(n)})}{\prod_{i=1}^{n-1} \left(1 - \frac {z_{\sigma(i)}}{z_{\sigma(i+1)}q_2} \right)} \prod^{\sigma(i) > \sigma(j)}_{1\leq i < j \leq n} \frac {\zeta \left(\frac {z_{\sigma(j)}}{z_{\sigma(i)}} \right)}{\zeta \left(\frac {z_{\sigma(i)}}{z_{\sigma(j)}} \right)} \prod_{i=1}^n Dz_i = 
$$
$$
\sum_{\sigma \in S(n)} \int_{|w_{\sigma^{-1}(1)}| \gg \dots \gg |w_{\sigma^{-1}(n)}|} \frac {r(w_{\sigma^{-1}(1)},\dots,w_{\sigma^{-1}(n)}) p(w_1^{-1}, \dots, w_n^{-1})}{\prod_{i=1}^{n-1} \left(1 - \frac {w_i}{w_{i+1}q_2} \right)} \prod^{\sigma(i) > \sigma(j)}_{1\leq i < j \leq n} \frac {\zeta \left(\frac {w_j}{w_i} \right)}{\zeta \left(\frac {w_i}{w_j} \right)}  \prod_{i=1}^n Dw_i
$$
where in the last equality we changed the variables to $w_i = z_{\sigma(i)}$. As we move the contours from $|w_{\sigma^{-1}(1)}| \gg \dots \gg |w_{\sigma^{-1}(n)}|$ toward $|w_1| \ll \dots \ll |w_n|$ in the integral above, we note that the only poles we might pick up are those of the form 
$$
\Big\{ w_i = w_{i+1}q_2 \Big\}_{i \in \{1,\dots,n-1\}}
$$
Thus, we conclude that the integral above is equal to
\begin{multline*}
\sum_{n_1+\dots+n_k = n} \int_{|x_1| \ll \dots \ll |x_k|} \left[ \underset{\{w_{n_1+\dots+n_{i-1}+1},\dots,w_{n_1+\dots+n_{i}}\}_{1\leq i \leq k} = \{x_iq_2^{n_i-1},\dots,x_i\}_{1\leq i \leq k}}{\text{Res}} \right. \\ \left. \sum_{\sigma \in S(n)}  \frac {r(w_{\sigma^{-1}(1)},\dots,w_{\sigma^{-1}(n)}) \cdot p(w_1^{-1}, \dots, w_n^{-1})}{\prod_{i=1}^{n-1} \left(1 - \frac {w_i}{w_{i+1}q_2} \right)} \prod^{\sigma(i) > \sigma(j)}_{1\leq i < j \leq n} \frac {\zeta \left(\frac {w_j}{w_i} \right)}{\zeta \left(\frac {w_i}{w_j} \right)} \right] \prod_{i=1}^k Dx_i
\end{multline*}
The second line of the expression above is
$$
\frac {p(w_1^{-1}, \dots, w_n^{-1})}{\prod_{i=1}^{n-1} \left(1 - \frac {w_i}{w_{i+1}q_2} \right) \prod_{1\leq i < j \leq n} \zeta \left(\frac {w_i}{w_j}\right)} \sum_{\sigma \in S(n)} r(w_{\sigma^{-1}(1)},\dots,w_{\sigma^{-1}(n)}) \prod_{\sigma(i) > \sigma(j)} \zeta \left(\frac {w_j}{w_i} \right)
$$
which directly implies \eqref{eqn:claim pairing} for $R$ of the form \eqref{eqn:rho}, as we needed to show.

\end{proof}

\end{proof}

\subsection{} In \cite{Shuf}, we introduced \textbf{slope} subalgebras
\begin{equation}
\label{eqn:slope}
\CB^\mu \subset \CS_{\loc}
\end{equation}
for any $\mu \in \BQ$, which are isomorphic to the algebra
$$
\Lambda = \BQ(q_1,q_2)[p_1,p_2,\dots]
$$
of symmetric polynomials in infinitely many variables (above, $p_t$ is interpreted as the $t$-th power sum function). Explicitly, we have an algebra isomorphism
\begin{equation}
\label{eqn:iso slope}
\tau^{\frac dn} : \Lambda \xrightarrow{\sim} \CB^{\frac dn}
\end{equation}
for any coprime $(n,d) \in \BN \times \BZ$, determined by the assignment
$$
\tau^{\frac dn}(p_t) = \oP_{nt,dt}
$$
If we let $h_t$ denote the $t$-th complete symmetric function, the power series identity
$$
1 + \sum_{t=1}^{\infty} \frac {h_t}{x^t} = \exp \left( \sum_{t=1}^{\infty} \frac {p_t}{tx^t} \right)
$$
and formula \eqref{eqn:oh to op} imply that $\tau^{\frac dn}(h_t) = \oH_{nt,dt}$ for all coprime $(n,d)$ and all $t \in \BN$.

\medskip

\begin{remark}

For any $\mu \in \BQ$, the isomorphism \eqref{eqn:iso slope} allows one to transport the usual Hall coproduct on $\Lambda$ to a coproduct $\Delta_\mu$ on $\CB^\mu$; the latter coproduct was given a shuffle algebra interpretation in \cite{Shuf}. In particular, this allows us to prove that
\begin{equation}
\label{eqn:f to h}
\oH_{n,0} = q_1^{\frac {n(n-1)}2} F_n 
\end{equation}
as both LHS and RHS are uniquely determined by the fact that they are group-like for $\Delta_0$ in $\CB^0$, and are annihilated by the linear map $\varphi$ (with $q_1 \leftrightarrow q_2$) of \loccit 

\end{remark}

\medskip

\subsection{} The ring $\Lambda$ is rich in automorphisms, as one can rescale the generators $p_t$ independently and arbitrarily. We will refer to the particularly important rescaling
$$
p_t' = p_t(q_1^t-1)
$$
as a \textbf{plethysm}. We therefore obtain elements $h_t' \in \Lambda$ via the usual formula
$$
1 + \sum_{t=1}^{\infty} \frac {h_t'}{x^t} = \exp \left( \sum_{t=1}^{\infty} \frac {p_t'}{tx^t} \right)
$$
It was shown in \cite[Section 2.3]{AGT} that 
$$
\tau^{\frac dn}(h_t') = \oH'_{nt,dt}
$$
for all coprime $(n,d) \in \BN \times \BZ$ and all $t \in \BN$, where
\begin{equation}
\label{eqn:hhnd}
\oH'_{n,d} = (q_1-1)^n(q_2-1)^n \cdot \Sym \left[\frac {\prod_{i=1}^n z_i^{\left \lfloor \frac {id}n \right \rfloor - \left \lfloor \frac {(i-1)d}n \right \rfloor}}{\prod_{i=1}^{n-1} \left(1-\frac {z_{i+1}}{z_iq_3}\right)} \prod_{1\leq i < j \leq n} \zeta \left(\frac {z_i}{z_j} \right) \right]
\end{equation}
\footnote{Note that our  $\oH'_{n,d}$ is $(1-q_1)H_{n,d}$ of \cite{GN}.} for any $(n,d) \in \BN \times \BZ$. Moreover, one can associate ribbon skew Schur functions
$$
s'_{\varepsilon} \in \Lambda
$$
to any sequence $\varepsilon$ consisting of zeroes and ones \footnote{Explicitly, for any $\varepsilon = (\varepsilon_1,\dots,\varepsilon_{t-1})$, we define $s'_{\varepsilon}$ as the skew Schur function associated to the size $t$ skew Young diagram whose first box is arbitrary, and whose $i+1$-th box is either to the right or below the $i$-th box, depending on whether $\varepsilon_i$ is 0 or 1 (see \cite[Section 6.12]{Shuf} for details).}, completely determined by
\begin{equation}
\label{eqn:rule}
s'_{\varepsilon} s'_{\varepsilon'} = s'_{\varepsilon 0 \varepsilon'} + s'_{\varepsilon 1 \varepsilon'}
\end{equation}
and the normalization $s'_{(0^{t-1})} = h_t'$ for all $t$. It was shown in \cite[Section 6.12]{Shuf} that
$$
\tau^{\frac dn}(s'_\varepsilon) = S'_{n,d, \varepsilon}
$$
where
\begin{multline}
\label{eqn:skew schur}
S'_{n,d, (\varepsilon_1 \dots \varepsilon_{t-1})} = (q_1-1)^n(q_2-1)^n \cdot \\ \Sym \left[\frac {\prod_{i=1}^n z_i^{\left \lfloor \frac {id}n \right \rfloor - \left \lfloor \frac {(i-1)d}n \right \rfloor} \prod^{1 \leq s \leq t-1}_{\varepsilon_s = 1} \left( - \frac {z_{as+1}}{z_{as}q_3} \right)}{\prod_{i=1}^{n-1} \left(1-\frac {z_{i+1}}{z_iq_3}\right)} \prod_{1\leq i < j \leq n} \zeta \left(\frac {z_i}{z_j} \right) \right]
\end{multline}
for any $(n,d) \in \BN \times \BZ$ with $\gcd(n,d) = t$ and $a = \frac nt$.

\medskip

\begin{remark} 

Comparing \eqref{eqn:skew schur} with \eqref{eqn:end 1} yields the identity
$$
\oH_{n,d} = \sum_{\varepsilon_1,\dots,\varepsilon_{t-1} \in \{0,1\}} \frac {q_1^{\sum^{1\leq s \leq t-1}_{\varepsilon_s = 0} s}}{(q_1-1) \dots (q_1^t-1)} \cdot S'_{n,d,(\varepsilon_1 \dots \varepsilon_{t-1})}
$$
which is simply $\tau^{\frac dn}$ applied to the symmetric function identity
$$
h_t = \sum_{\varepsilon_1,\dots,\varepsilon_{t-1} \in \{0,1\}} \frac {q_1^{\sum^{1\leq s \leq t-1}_{\varepsilon_s = 0} s}}{(q_1-1) \dots (q_1^t-1)} \cdot s'_{(\varepsilon_1 \dots \varepsilon_{t-1})}
$$

\end{remark}

\noindent The upshot of the discussion above is that the elements $\{\oH_{nt,dt}, \oH'_{nt,dt}\}_{t \in \BN}$ play the roles of ones and the same symmetric functions for any coprime $(n,d) \in \BN \times \BZ$, under the isomorphisms \eqref{eqn:iso slope}. Thus, understanding these elements for one slope (say 0) would yield an understanding for all slopes.

\end{document}